\documentclass[12pt]{amsart}
\usepackage{latexsym}
\usepackage{amsmath}
\usepackage{amsfonts}
\usepackage{amssymb}
\usepackage{amsthm}

\theoremstyle{plain}
\newtheorem{theorem}{Theorem}[section]
\newtheorem{lemma}[theorem]{Lemma}
\newtheorem{proposition}[theorem]{Proposition}

\newtheorem*{theorem*}{Theorem}

\newtheorem{corollary}[theorem]{Corollary}

\theoremstyle{definition}
\newtheorem{definition}[theorem]{Definition}
\newtheorem*{example}{Example}

\theoremstyle{remark}

\font\tenmsb=msbm10 at 12pt \font\sevenmsb=msbm7 at 8pt
\font\fivemsb=msbm5 at 6pt
\newfam\msbfam
\textfont\msbfam=\tenmsb \scriptfont\msbfam=\sevenmsb
\scriptscriptfont\msbfam=\fivemsb

\newcommand{\bbE}{\mathbb{E}}
\newcommand{\bbN}{\mathbb{N}}

\newcommand{\bbZ}{\mathbb{Z}}
\newcommand{\norm}[1]{\lvert\!| #1~\!|\!|}
\newcommand{\nnorm}[1]{\lvert\!|\!| #1\!|\!|\!|}

\begin{document}

\title[Convergence of Several Commuting Polynomial Averages]{Convergence of
Polynomial Ergodic Averages of Several Variables for Some Commuting Transformations}
\author{Michael C. R. Johnson}
\maketitle
  {\footnotesize 
  \centerline{ Department of Mathematics }
  \centerline{ Northwestern  University}
  \centerline{ Evanston, IL 60201} }

\begin{abstract}
Let $(X,\mathcal{B},\mu)$  be a probability space and let $T_1,\ldots , T_l$ be $l$
commuting invertible measure preserving transformations \linebreak of $X$. We show
that if $T_1^{c_1} \ldots T_l^{c_l}$ is ergodic for each $(c_1,\ldots,c_l)\neq (0,\ldots,0)$, then the averages $\frac{1}{|\Phi_N|}\sum_{u\in\Phi_N}\prod_{i=1}^r T_1^{p_{i1}(u)}\ldots T_l^{p_{il}(u)}f_i$ converge in $L^2(\mu)$ for all polynomials $p_{ij}\colon \mathbb{Z}^d\to\mathbb{Z}$,
all $f_i\in L^{\infty}(\mu)$, and all F\o lner sequences $\{\Phi_N\}_{N=1}^{\infty}$ in
$\mathbb{Z}^d$.
\end{abstract}
\medskip
\section{introduction}

In 1996, Bergelson and Leibman proved the following generalization of Furstenberg's
Multiple Recurrence Theorem \cite{F77}, corresponding to the multidimensional
polynomial version of Szemer\'edi's theorem.
\begin{theorem}\cite{BL}\label{X}
Let $(X,\mathcal{B},\mu)$ be a probability space, let $T_1, \ldots, T_l$ be
commuting invertible measure preserving transformations of $X$, let $p_{ij}\colon
\bbZ \to \bbZ$ be polynomials satisfying $p_{ij}(0)=0$ for all $1\leq i\leq r, 1\leq
j\leq l$, and let $A\in \mathcal{B}$ with $\mu(A)>0$. Then
$$ \liminf_{N\to\infty}\frac{1}{N}\sum_{n=0}^{N-1}\mu(\bigcap_{i=1}^r
T_1^{-p_{i1}(n)}\ldots T_l^{-p_{il}(n)}A)>0.$$
\end{theorem}
Furstenberg's theorem corresponds to the case that $p_{ij}(n)=n$ for $i=j$,
$p_{ij}(n)=0$ for $i\neq j$ and each $T_i=T_1^i$. In this linear case, Host and Kra
\cite{HK02} showed that the lim inf is in fact a limit. Host and Kra \cite{HK04} and
Leibman \cite{Le4} proved convergence in the polynomial case assuming all $T_i=T_1$.
It is natural to ask whether the general commuting averages for polynomials in Theorem \ref{X}
converge.

\begin{definition} We say $(T_1, \ldots, T_l)$ is a {\bf totally ergodic generating set} of  invertible measure preserving transformations of $X$ if each $T_1^{c_1}T_2^{c_2}\dots T_l^{c_l}$ is ergodic for any choice of $(c_1,\ldots,c_l)\neq (0,\ldots,0)$. 
\end{definition}
We note that if $(T_1, \ldots, T_l)$ is a totally ergodic generating set of  invertible measure preserving transformations of a non-trivial probability space $(X,\mathcal{B},\mu)$, then the associated group of transformations  generated by $T_1, \ldots, T_l$ is a free abelian group with $l$ generators.
We show that given a totally ergodic generating set of transformations, we obtain convergence in $L^2(\mu)$ for the averages in Theorem \ref{X}. We prove a statement replacing indicator functions with arbitrary functions in $L^{\infty}(\mu)$.

\begin{theorem}\label{A}
Let $(X,\mathcal{B},\mu)$ be a probability space, let $(T_1, \ldots, T_l)$ be a
totally ergodic generating set of commuting invertible measure preserving transformations of
$X$, and let $p_{ij}\colon \bbZ^d \to \bbZ$ for $1\leq i\leq r, 1\leq j\leq l$ be
polynomials. For any $f_1,\ldots,f_r \in L^{\infty}(\mu)$ and any F\o lner sequence
$\{\Phi_N\}_{N=1}^{\infty}$ in $\bbZ^d$, the averages 
\begin{multline}
\label{1} \frac{1}{|\Phi_N|}\sum_{u\in \Phi_N}\prod_{i=1}^r
f_i(T_1^{p_{i1}(u)}\ldots T_l^{p_{il}(u)}x)
\end{multline}
converge in $L^2(\mu)$ as $N\to\infty$.
\end{theorem}

Without the assumption that $(T_1,\ldots, T_l)$ form a totally ergodic generating set, convergence for the above averages in \eqref{1} remains open and is only known in
the linear case. Frantzikinakis and Kra \cite{FrK} showed that
given $p_{ij}(n)=n$ for $i=j$ and $p_{ij}(n)=0$ for $i\neq j$, if we assume that
$T_i$ is ergodic for each $i\in\{1,\ldots,l\}$ and $T_iT_j^{-1}$ is ergodic for all
$i\neq j$, we obtain convergence in $L^2(\mu)$. Tao \cite{Ta} recently proved convergence in $L^2(\mu)$ for the general linear case without the ergodicity assumptions needed in \cite{FrK}.

In previous results, convergence was often shown by proving that the averages in (1) do
not change by replacing each function with its conditional expectation on a certain
characteristic factor, namely an inverse limit of nilsystems. This characteristic
factor, is then shown to have algebraic structures for which convergence is known.
We define these terms precisely in the section below. To prove our theorem, we
combine this technique with  a modified version of PET-induction as introduced by Bergelson \cite{Be}.

\section{Preliminaries}

For simplicity, we assume all functions are real valued. All theorems and
definitions hold for complex valued functions with obvious minor modifications.
Throughout, we use the notation $Tf=f(T)$.
\subsection{Nilsystems}

\begin{definition}
Let $G$ be a $k$-step nilpotent Lie group, let $\Gamma$ be a discrete cocompact
subgroup of $G$, let $X=G/\Gamma$, and let $\mathcal{B}$ be the Borel $\sigma$-algebra associated to $X$. For each $g\in G$, let $T_g :G/\Gamma \to
G/\Gamma$ be defined by $T_g(x\Gamma)=gx\Gamma$, and let $\mu$ be Haar measure, the
unique measure on $(X,\mathcal{B})$ invariant under left translations by elements in G. We call $(X,\mathcal
B,\mu, (T_g,g\in G))$ a \bf{nilsystem}.
\end{definition}

\begin{definition}
A sequence of finite subsets $\{\Phi_N\}_{N=1}^{\infty}$ of a countable, discrete group G is a {\bf F\o
lner sequence} if for all $g\in G$, 
$$\lim_{n\to\infty} \frac{|g \Phi_n \Delta \Phi_n|}{|\Phi_n|} = 0,$$ 
where $\Delta$ is the symmetric difference operation.
\end{definition}
Ergodic averages in nilsystems have been well studied. We make use of the
following theorem of Leibman:

\begin{theorem}{\cite{Le1}}\label{B}
Let $(X,\mathcal B,\mu,(T_g, g\in G))$ be a nilsystem with $X=G/\Gamma$, $g_{_1},\ldots,g_{_l}
\in G$, and $p_1,\ldots , p_l\colon \bbZ^d \to \bbZ$ be polynomials. Then for any
$f\in C(X)$ and any F\o lner sequence $\{\Phi_N\}_{N=1}^{\infty}$ in $\bbZ^d$, the
averages 

$$\frac{1}{|\Phi_N|}\sum_{u\in \Phi_N}T_{g_{_1}}^{p_1(u)}\ldots T_{g_{_l}}^{p_l(u)}f$$ 
converge pointwise as $N\to\infty$.
\end{theorem}
\begin{corollary}\label{C}
Let $(X,\mathcal B,\mu,(T_g,g\in G))$ be a nilsystem with $X=G/\Gamma$,
$g_{_1},\ldots,g_{_l} \in G$, and $p_{ij}\colon\bbZ^d \to \bbZ$ for $1\leq i\leq r,
1\leq j\leq l$ be polynomials. Then for any $f_1,\ldots,f_r \in L^{\infty}(\mu)$ and
any F\o lner sequence $\{\Phi_N\}_{N=1}^{\infty}$ in $\bbZ^d$, the averages 
$$\frac{1}{|\Phi_N|}\sum_{u\in \Phi_N}\prod_{i=1}^r T_{g_{_1}}^{p_{i1}(u)}\ldots
T_{g_{_l}}^{p_{il}(u)}f_i$$ 
converge in $L^2(\mu)$ as $N\to\infty$. 
\end{corollary}
\begin{proof}
We apply Theorem \ref{B} to $X^{r}$, with transformations $\hat{T}_{ij}\colon X^r\to X^r$ for $1\leq i\leq r,
1\leq j\leq l$ defined by
$$\hat{T}_{ij}(x_1,x_2,\ldots, x_r)=(x_1,\ldots,x_{i-1}, T_{g_{_j}}(x_i),x_{i+1},\ldots,x_r).$$ Using polynomials $p_{ij}\colon\bbZ^d \to \bbZ$ for $1\leq i\leq r,
1\leq j\leq l$ and
 $f=f_1\otimes \ldots \otimes f_r$, we get 
 $$\hat{T}_{11}^{p_{11}(u)}\ldots \hat{T}_{1l}^{p_{1l}(u)}\ldots \hat{T}_{r1}^{p_{r1}(u)}\ldots \hat{T}_{rl}^{p_{rl}(u)}f=\prod_{i=1}^r T_{g_{_1}}^{p_{i1}(u)}\ldots
T_{g_{_l}}^{p_{il}(u)}f_i.$$

Theorem \ref{B} guarantees the required averages converge pointwise for each $f_1,\ldots, f_r\in
C(X)$. Using the density of $C(X)$ in $L^{2}(\mu)$, $L^2(\mu)$ convergence follows for
arbitrary $f_1,\ldots, f_r\in L^{\infty}(\mu)$.
\end{proof}
\subsection{The Host-Kra seminorms $\nnorm {\cdot}_k$} We briefly review the
construction of the Host-Kra seminorms on $L^{\infty}(\mu)$ from \cite{HK02}. As our
setting deals with multiple commuting transformations, we must specify which
transformation is used. In this section, $T$ is an ergodic measure preserving transformation of $(X,\mathcal{B},\mu)$.

For each $k\geq 0$ we define a probability measure $\mu^{[k]}_T$ on $X^{[k]} =
X^{2^k}$, invariant under $T^{[k]}= T\times \dots \times T$ ($2^k$ times). 

Set $\mu^{[0]}_T=\mu$. For $k\geq 0$, let $\mathcal{I}^{[k]}_T$ be the
$\sigma$-algebra of $T^{[k]}$-invariant subsets of $X^{[k]}$. Then define
$\mu^{[k+1]}_T= \mu^{[k]}_T \times_{\mathcal{I}_T^{[k]}}\mu^{[k]}_T$ to be the
relatively independent square of $\mu^{[k]}_T$ over $\mathcal{I}_T^{[k]}$. This
means for $F,G\in L^{\infty}(\mu^{[k]})$

$$\int_{X^{[k+1]}}F({\bf x'})G({\bf x''}) d\mu_T^{[k+1]}({\bf x'},{\bf x''}) =
\int_{X^{[k]}}\bbE (F|\mathcal{I}_T^{[k]}) \bbE (G|\mathcal{I}_T^{[k]})
d\mu^{[k]}_T,$$ where $\mathbb{E}(\cdot | \cdot)$ is the conditional expectation operation.

Using these measures, define 
$$\nnorm f_{k,T}^{2^k} = \int_{X^{[k]}} \prod_{j=0}^{2^k-1}f(x_j) d\mu^{[k]}_T({\bf
x})$$ 
for a bounded function $f \in L^{\infty}(\mu)$ and $k\geq 1$. It is shown in
\cite{HK02} that for every $k\geq 1$ and every ergodic $T$, $\nnorm {\cdot}_{k,T}$
is a seminorm on $L^{\infty}(\mu)$. Also, for $f\in L^{\infty}(\mu)$, we have
$\nnorm f_{1,T} = |\int f d\mu |$ and for every $k \geq 1$, $\nnorm f_{k,T} \leq
\nnorm f_{k+1,T} \leq \norm f_{L^{\infty}(\mu)}.$

\subsection{The Host-Kra factors $Z_k(X)$} We now define an increasing sequence of
factors $\{Z_k(X,T) : k\geq 0\}$ as constructed in \cite{HK02}. Let
$\mathcal{Z}_k(X,T)$ be the $T$-invariant sub-$\sigma$-algebra characterized by the
following property: for every $f\in L^{\infty}(\mu)$, $\bbE (f |\mathcal{Z}_k(X,T))
= 0$ if and only if $\nnorm f_{k+1,T} = 0$. We define $Z_k(X,T)$ to be the factor of
$X$ associated to the sub-$\sigma$-algebra $\mathcal{Z}_k$. Thus $Z_0(X,T)$ is the
trivial factor and $Z_1(X,T)$ is the Kronecker factor. {\it A priori}, these
constructions depend on the transformation $T$.

Indeed, the following observation of Frantzikinakis and Kra shows that given basic
assumptions, none of the previous constructions depend on the transformation $T$. 
\begin{proposition}\cite{FrK}\label{D}
Assume that T and S are ergodic commuting invertible measure preserving
transformations of a space $(X,\mathcal{B},\mu)$. Then for all $k\geq 1$ and all
$f\in L^{\infty}(\mu)$, $\nnorm f_{k,T} = \nnorm f_{k,S}$ and $Z_k(X,T) = Z_k(X,S)$.
\end{proposition}
Thus we discard $T$ from our notation.

\begin{definition}
We call a probability space $(X,\mathcal{B},\mu)$ with $l$ invertible commuting
measure preserving transformations $T_1,\ldots,T_l$, an {\bf (invertible commuting
measure preserving) system}. If $(T_1,\ldots,T_l)$ is also a totally
ergodic generating set, then we call it a {\bf  freely generated totally ergodic system (with generators $(T_1,\ldots,T_l)$)}. We denote it as
$(X,\mathcal{B},\mu,(T_1,\ldots,T_l))$. A system $(X,\mathcal{B},\mu,(T_1,\ldots,T_l))$
is an {\bf inverse limit} of systems $(X,\mathcal{B}_i,\mu_i,(T_1,\ldots,T_l))$ if
each $\mathcal{B}_i\subset\mathcal{B}_{i+1}$ and $\mathcal{B}=\bigvee_{i=1}^{\infty}
\mathcal{B}_i$ up to sets of measure zero.
\end{definition}

The main result of the Host-Kra theory is that each of the factors $(Z_k,T_i)$ is
isomorphic to an inverse limit of $k$-step nilsystems. However, such isomorphism {\it
a priori} depends on the transformation $T_i$. (Note that by Proposition \ref{D}, $Z_k(X,T_i)$,
does not depend on $i$). In \cite{FrK}, they deal specifically
with this technicality.  We say that a system $(X,\mathcal{B},\mu,(T_1,\ldots,T_l))$
has {\bf order} $k$ if $X=Z_k(X)$. 

\begin{theorem}{\cite{FrK}}\label{E}
Any system $(X,\mathcal{B},\mu,(T_1,\ldots,T_l))$ of order k is an inverse limit of a
sequence of systems $(X,\mathcal{B}_i,\mu_i,(T_1,\ldots,T_l))$, each arising from $k$-step nilsystems, where $X=G_i/\Gamma_i$ and each transformation $T_1,\ldots,T_l$ is a left translation of $G_i/\Gamma_i$ by an element in $G_i$.
\end{theorem}
By combining Theorem \ref{E} and Corollary \ref{C}, Theorem \ref{A} is proved in the
case that $X=Z_k(X)$ for some $k$. 

\subsection{Characteristic factors and ED-sets}
\begin{definition}
We say a sub-$\sigma$-algebra $\mathcal{X}\subseteq\mathcal{B}$ is a {\bf
characteristic factor for $L^2(\mu)$-convergence} of the averages 
$$\hskip -2.3truein(\ref{1}) \hskip .3truein \frac{1}{|\Phi_N|}\sum_{u\in \Phi_N}\prod_{i=1}^r T_1^{p_{i1}(u)}\ldots
T_l^{p_{il}(u)}f_i$$
if $\mathcal{X}$ is $T_j$ invariant for all $1\leq j \leq l$ and the averages in
(\ref{1}) converge to 0 in $L^2(\mu)$ for any F\o lner sequence $\{\Phi_N\}_{N=1}^{\infty}$ in $\bbZ^d$ whenever $\bbE(f_i|\mathcal{X})=0$ for some
$1\leq i \leq r$. 
\end{definition} 

Using the multilinearity of our averages in (1), it only remains to show that for
some $k\in\bbN$, $Z_k(X)$ is a characteristic factor. 

To simplify future arguments, we require that our set of polynomials have a property related to being essentially distinct, as defined in \cite{Le4}.
\begin{definition}
We say the set of polynomials $P=\{p_{ij}\colon\bbZ^d \to \bbZ$ for $1 \leq i \leq r, 1
\leq j\leq l\}$ is an {\bf ED-set} if all of the following hold:
\begin{enumerate}
\item Each $p_{ij}$ in $P$ is not equal to a nonzero constant.
\item No two polynomials $p_{i_1j_1},p_{i_2j_2}$ in $P$ differ by a nonzero constant.
\item For each $i=1,\ldots,r$, there is some $j\in \{1,\ldots, l\}$ where $p_{ij}$ is nonzero.
\item For each distinct pair $i_1,i_2\in\{1,\ldots,r\}$, there is some $j\in \{1,\ldots, l\}$ where $p_{i_1j}\neq p_{i_2j}$.
\end{enumerate} Conditions (1) and (2) are related to the polynomials being essentially distinct. When $P$ is viewed as an $r\times l $ matrix whose entries are polynomials, condition (3) requires that $P$ contains no rows of all zeros, and condition (4) requires that $P$ does not have identical rows.
\end{definition}

We note that Theorem \ref{A} is trivially true if all the polynomials are identically zero.
By replacing each $f_i$ with $T_1^{c_1}\ldots T_l^{c_l}f_i$ for
some $c_1,\ldots,c_l\in\bbZ$, we may assume that our set of polynomials satisfies conditions (1) and (2). When $T_1^{p_{i1}}\ldots T_l^{p_{il}}f_i=f_i$, we  factor $f_i$ out of our average. Thus, we further  assume our polynomials satisfy condition (3). By writing $T_1\ldots T_l f T_1\ldots T_l g$ as $T_1\ldots T_l (fg)$ we may assume that our set of polynomials also satisfies condition (4), and hence is an ED-set.
 Thus the main theorem is a consequence of the following:

\begin{proposition}\label{R}
Let $(X,\mathcal{B},\mu,(T_1,\ldots,T_l))$ be a freely generated totally ergodic system and
$P=\{p_{ij}\colon\bbZ^d \to \bbZ$ for $1 \leq i \leq r, 1 \leq j\leq l\}$ be an ED-set of polynomials. Then there exists $k\in\bbN$ such that for any \linebreak $f_1, \ldots , f_r \in L^{\infty}(\mu)$ with
$\nnorm {f_m}_k = 0$ for some $1\leq m \leq r$, we have $$\limsup_{N\to\infty}
\Big\|\frac{1}{|\Phi_N|}\sum_{u \in \Phi_N} (\prod_{i=1}^r
T_1^{p_{i1}(u)}T_2^{p_{i2}(u)}\ldots T_l^{p_{il}}f_i )\Big\|_{L^2(\mu)}=0$$ for any
F\o lner sequence $\{\Phi_N\}_{N=1}^{\infty}$ in $\bbZ^d$.
\end{proposition}

We note that the above integer $k$ is only dependent on the set of polynomials
$P$ and not on the system $(X,\mathcal{B},\mu, (T_1,\dots,T_l))$ or the
dimension $d$. By relabeling our polynomials and functions, we need only prove
Proposition \ref{R} in the case that $\nnorm{f_1}_k=0$ for some $k\in\bbN$.

\section{Linear case}

To prove proposition \ref{R}, we use PET-induction as introduced by Bergelson in
\cite{Be}. In this section we prove the base case of the induction. 

\begin{proposition}\label{G}
Let $(X,\mathcal{B},\mu,(T_1,\ldots,T_l))$ be a freely generated totally ergodic system and
$P=\{p_{ij}\colon\bbZ^d \to \bbZ$ for $1\leq i\leq r, 1\leq j\leq l\}$ be an ED-set of linear functions. Then there exists a constant $C>0$ dependent only on the
set of polynomials, such that
$$\limsup_{N\to\infty} \Big\|\frac{1}{|\Phi_N|}\sum_{u \in \Phi_N} (\prod_{i=1}^r
T_1^{p_{i1}(u)}T_2^{p_{i2}(u)}\ldots T_l^{p_{il}(u)}f_i )\Big\|_{L^2(\mu)}$$
$$\leq C \min_{1\leq i\leq r} \nnorm {f_i~}_{r+1}$$
for any $f_1,\ldots ,f_r \in L^{\infty}(\mu)$ with $\|f_i\|_{L^{\infty}(\mu)}\leq 1$ and any F\o lner
sequence $\{\Phi_N\}_{N=1}^{\infty}$ in $\bbZ^d$.
\end{proposition}

As a corollary, we get that $Z_r(X)$ is characteristic for the averages in
(\ref{1}) when each of the polynomials in $P$ is linear. We use the following version of the van der Corput lemma in the inductive
process to reduce each average to a previous step. 

\begin{lemma}{\cite{BMZ}}\label{F}
Let $\{g_u\}_{u\in G}$ be a bounded family of elements of a Hilbert space
$\mathcal{H}$ indexed by elements of a finitely generated abelian group G and let
$\{\Phi_N\}_{N=1}^{\infty}$ be a F\o lner sequence in G. 
\begin{enumerate}
\item For any finite set $F\subseteq G$,
$$\limsup_{N\to\infty}\Big\|\frac{1}{|\Phi_N|}\sum_{u\in \Phi_N} g_u\Big\|^2 \leq
\limsup_{N\to\infty}\frac{1}{|F|^2}\sum_{v,w\in F}\frac{1}{|\Phi_N|}\sum_{u\in
\Phi_N}\langle g_{u+v},g_{u+w} \rangle.$$

\item There exists a F\o lner sequence $\{\Theta_M\}_{M=1}^{\infty}$ in $G^3$ such
that 
$$\limsup_{N\to\infty}\Big\|\frac{1}{|\Phi_N|}\sum_{u\in \Phi_N} g_u\Big\|^2 \leq
\limsup_{M\to\infty}\frac{1}{\Theta_M}\sum_{(u,v,w)\in \Theta_M}\langle
g_{u+v},g_{u+w} \rangle.$$
\end{enumerate}
\end{lemma}
Leibman proved the following lemma in his proof of convergence for a single
transformation \cite{Le4}. We likewise use his lemma to prove the linear case for
multiple commuting transformations.

\begin{lemma}\cite{Le4}\label{H}
\begin{enumerate}
\item Let $p_i\colon\bbZ^d \to \bbZ$ be nonconstant linear functions for each
$i=1,\ldots,l$. There exists a constant $C$, such that for any $f \in
L^{\infty}(\mu)$ and any F\o lner sequence $\{\Phi_N\}_{N=1}^{\infty}$ in $\bbZ^d$, 
$$\lim_{N\to\infty} \Big\|\frac{1}{|\Phi_N|}\sum_{u\in \Phi_N}T_1^{p_{1}(u)}\ldots
T_l^{p_{l}(u)}f\Big\|_{L^2(\mu)} \leq C\nnorm f _2.$$
\item Let $p_i\colon\bbZ^d \to \bbZ$ be nonconstant linear functions for each
$i=1,\ldots,l$. There exists a constant C, such that for any $f \in L^{\infty}(\mu)$
and any F\o lner sequence $\{\Phi_N\}_{N=1}^{\infty}$ in $\bbZ^d$,
$$\lim_{N\to\infty} \frac{1}{|\Phi_N|} \sum_{u\in \Phi_N} \nnorm{f \cdot
T_1^{p_{1}(u)}\ldots T_l^{p_{l}(u)}f}_k^{2^k} \leq C \nnorm f_{k+1}^{2^{k+1}}.$$
\end{enumerate}
\end{lemma}
We note here that  Lemma \ref{H} is similar to Lemmas 7 and 8 in \cite{Le4} but
with multiple commuting transformations. The only step needed to alter his proof is
to show our average also convergences to the conditional expection of $f$ onto the
appropriate sub-$\sigma$-algebra. But this follows from classical results on
convergence for amenable group actions.

\begin{proof}[Proof of Proposition~\ref{G}]
To simplify notation, we write $T_1^{p_{i1}(u)}\ldots T_l^{p_{il}(u)}$ as
$S^{p_i(u)}$. Since each $p_{ij}$ is a linear polynomial, we have
$S^{p_i(u)}S^{p_i(v)}=S^{p_i(u+v)}$.

We proceed by induction on $r$. For $r=1$, we are done by Lemma \ref{H}. Assume the
proposition holds for $r-1$ functions. Let $f_1, \ldots, f_r$ be essentially bounded
functions on $X$ with $\norm{f_i}_{L^{\infty}(\mu)} \leq 1$ for all $1 \leq i \leq
r$, and let $\{\Phi_N\}_{N=1}^{\infty}$ be a F\o lner sequence in $\bbZ^d$. By
applying Lemma \ref{F} to $g_u=S^{p_i(u)}f_1 \ldots S^{p_r(u)}f_r$, for any finite
$F\subseteq \bbZ^d$, we get
\begin{multline*}
\limsup_{N\to\infty} \Big\|\frac{1}{|\Phi_N|}\sum_{u\in \Phi_N}\prod_{i=1}^r
S^{p_i(u)}f_i\Big\|_{L^2(\mu)}^2\\
\leq \limsup_{N\to\infty} \frac{1}{|F|^2}\sum_{v,w\in F}
\frac{1}{|\Phi_N|}\sum_{u\in \Phi_N}\int_X\prod_{i=1}^r S^{p_i(u+v)}f_i \cdot \prod
_{i=1}^r S^{p_i(u+w)}f_i d\mu\\
= \limsup_{N\to\infty} \frac{1}{|F|^2}\sum_{v,w\in F} \frac{1}{|\Phi_N|}\sum_{u\in
\Phi_N}\int_X\prod_{i=1}^{r-1}S^{p_i(u)}S^{-p_r(u)}(S^{p_i(v)}f_i \\ 
\cdot S^{p_i(w)}f_i) \cdot(S^{p_r(v)}f_r \cdot S^{p_r(w)}f_r) d\mu\\
\leq \frac{1}{|F|^2}\sum_{v,w\in
F}\limsup_{N\to\infty}\Big\|\frac{1}{\Phi_N}\sum_{u\in
\Phi_N}\prod_{i=1}^{r-1}S^{(p_i-p_r)(u)}(S^{p_i(v)}f_i \cdot
S^{p_i(w)}f_i)\Big\|_{L^2(\mu)}.
\end{multline*}
Since $P$ is an ED-set, so is the family $\{(p_{ij}-p_{rj})\colon \bbZ^d\to\bbZ \text{ for }1\leq i \leq r-1,1\leq j \leq l\}$.
By the induction process, there exists a constant $C$, independent of $f_1, \ldots,
f_r$ and $\{\Phi_N\}_{N=1}^{\infty}$, such that 
\begin{multline*}\limsup_{N\to\infty}\Big\|\frac{1}{\Phi_N}\sum_{u\in
\Phi_N}\prod_{i=1}^{r-1}S^{(p_i-p_r)(u)}(S^{p_i(v)}f_i \cdot
S^{p_i(w)}f_i)\Big\|_{L^2(\mu)} \\
\leq C\nnorm{(S^{p_i(v)}f_i \cdot S^{p_i(w)}f_i)~}_r
\end{multline*}
for all $(v,w) \in \bbZ^{2d}$ and $i\in\{1,\ldots,r\}$. Thus for any finite set $F
\subset \bbZ^d$ and $i\in\{1,\ldots,r\}$, 
\begin{multline*}
\limsup_{N\to\infty} \Big\|\frac{1}{|\Phi_N|}\sum_{u\in \Phi_N}\prod_{i=1}^r
S^{p_i(u)}f_i\Big\|_{L^2(\mu)}\\
\leq \Bigl(\frac{C}{|F|^2}\sum_{v,w \in F}\nnorm{(S^{p_i(v)}f_i \cdot
S^{p_i(w)}f_i)~}_r\Bigr)^{1/2}\\
\leq C^{1/2}\Bigl(\frac{1}{|F|^2}\sum_{v,w \in F}\nnorm{(f_i \cdot
S^{p_i(w-v)}f_i)~}_r^{2^r}\Bigr)^{(1/2)^{r+1}}.
\end{multline*}

Let $\{\Psi_N\}_{N=1}^{\infty}$ be any F\o lner sequence in $\bbZ^d$. Thus
$\{\Psi_N\times \Psi_N\}_{N=1}^{\infty}$ is a F\o lner sequence in $\bbZ^{2d}$. By
Lemma \ref{H} we have for each $i\in\{1,\ldots,r\}$
$$\limsup_{M\to\infty}\frac{1}{|\Psi_M|^2}\sum_{v,w\in\Psi_M}\nnorm{f_i \cdot
S^{p_i(w-v)}f_i~}_r^{2^r} \leq c\nnorm{f_i~}_{r+1}^{2^{r+1}}$$ with c independent
of $f_i$. By replacing $F$ with $\Psi_N$ for each $N\in\bbN$, we get
$$\limsup_{N\to\infty} \Big\|\frac{1}{|\Phi_N|}\sum_{u\in \Phi_N}\prod_{i=1}^r
S^{p_i(u)}f_i\Big\|_{L^2(\mu)} \leq C^{1/2}c^{(1/2)^{r+1}}\min_{i\leq
r}\nnorm{f_i~}_{r+1}.$$

\end{proof}

\section{Non-linear Case} 

We now deal with the inductive step. A set of polynomials $P = \{p_{ij} : 1 \leq i
\leq r, 1 \leq j\leq l\}$ where each $p_{ij}\colon\bbZ^d \to \bbZ$ is called a {\bf
(integer) polynomial family}. We view $P$ as an $r\times l$ matrix whose entries are the polynomials $p_{ij}$. We define the {\bf degree} of a family $P$,
$$\deg(P)=\max \{ \deg(p_{ij})\colon p_{ij} \in P\}.$$ 
Let $D\in\bbN$.
We define the {\bf column degree}
of a polynomial family $P$ with $\deg(P)\leq D$ to be the vector $C(P)= (c_1,\ldots,c_D)$ where $c_i$ is the number of columns whose maximal degree is $i$.

We say that two polynomials $p$, $q$ are
equivalent if $\deg(p)=\deg(q)$ and $\deg(p-q) < \deg(p)$. Thus any collection of polynomials can be partitioned into equivalence classes.  We define the degree of an equivalence class of polynomials to be equal to the degree of any of its representatives. 

For a family $P$ with $\deg(P)\leq D$, we define the {\bf column weight} of a column $j$, to be the vector
$w_j(P)=(w_{1j},\ldots,w_{Dj})$ with each $w_{ij}$ equal to the number of
equivalence classes in $P$ of degree $i$ in column $j$. Given \linebreak two vectors ${\bf
v}=(v_1,\ldots,v_{D})$, ${\bf v'}=(v'_1\ldots, v'_{D})$, we say ${\bf v}<{\bf v'}$
 there exists $n_0$ such that $v_{n_0}<v'_{n_0}$ and for
each $n>n_0$, $v_n=v'_n$. 
Thus the set of weights and the set of column degrees
become well ordered sets. 

  For each polynomial family $P$ with $\deg(P)\leq D$, we define the {\bf subweight} of $P$ to be the matrix $w(P)=[w_1(P)\ldots w_D(P)]$ whose columns are the corresponding column weights of $P$. 
Due the fact that our polynomial family may have many polynomial entries that are zero, we must modify the PET-induction scheme from that of \cite{Le4}. We introduce the following notation to record the position of such zeros in $P$. Let

$I_0=\{i\in\{1,\ldots,r\} : p_{ij} = 0$ for all $j =1,\ldots, l\}$, 

$I_1=\{i\in\{1,\ldots,r\} : \deg(p_{ij}) \leq 1$ for all $j=1\ldots, l\}\setminus I_0$, and 

$I_2 =\{1,\ldots,r\} \setminus (I_0 \cup I_1)$. 

 When $P$ is an ED-set, $I_0$ is empty, $I_1$ records which nonzero rows contain only polynomials with degree at most 1, while $I_2$ records which rows contain a polynomial of degree greater than 2.
Define $H_0(P)=I_1\cup I_2$ and inductively define $$H_j(P)=\{i\in\{1,\ldots, r\}\colon
p_{ij}=0\}\cap H_{j-1}(P)$$ for $1\leq j\leq l-1$ (we omit the polynomial
family $P$ when there is no confusion which family we are dealing with). Thus, $H_j$
records which non-identically zero rows have zeros in columns
$1,\ldots,j$. Pick $j_0$ to be the smallest $j\geq 1$ such that
$H_{j}=\emptyset$. In the case that column $1$ has no zero entries, we
note that $j_0=1$.

For each polynomial family $P$ and integer $a=1,\ldots,l$, we define the
sub-polynomial family $$P^{a}=\{p_{ij}\colon i\in H_{a-1}(P),a\leq j\leq l\}.$$
We note that the entries in the first column in $P^{a}$ are precisely the entries of column $a$ of $P$ from nonzero rows whose polynomials are all identically zero in columns $1,\ldots,a-1$. We note that when $P$ is an ED-set, $P^1=P$.

For each polynomial family $P$ with $\deg(P)\leq D$, we define the {\bf weight} of $P$ to be the ordered set of matrices $W(P)=\{w(P^1),\ldots, w(P^l)\}$. Given two polynomial families $P$ and $Q$ where $\deg(P)$, $\deg(Q)\leq D$, we say that $W(Q)<W(P)$ if there exists $J,A\in\{1,\ldots, l\}$ such that $w_J(Q^A)<w_J(P^A)$, but $w_J(Q^a)=w_J(P^a)$ for all $1\leq a<A$ and $w_j(Q^a)=w_j(P^a)$ for all $1\leq j<J$ and $a=1,\ldots l$.

\begin{example}\label{Example}
Let $P=\left(\begin{matrix}
n^2&2n&n\\
0&n^2&0\\
0&2n^2&3n
\end{matrix}\right)$. We see that $P$ is an ED-set,  and $H_1(P)=\{2,3\}$. Thus $P^2=\left(\begin{matrix}
n^2&0\\
2n^2&3n
\end{matrix}\right)$. Since $H_2(P)=\emptyset$, $P^3$ is the empty family. Therefore $w(P^1)=\left[\begin{matrix} 0&1&2\\
1&2&0\end{matrix}\right]$, $w(P^2)=\left[\begin{matrix} 0&1\\
2&0\end{matrix}\right]$, and $w(P^3)=\left[\begin{matrix} 0\\
0\end{matrix}\right]$.
Let $$Q=\left(\begin{matrix}
n^2-2n+1&-n^2+1&n+1\\
n^2+2n+1&-n^2+1&n+1\\
0&-4n&0\\
0&n^2-6n+1&3n+3\\
0&n^2+2n+1&3n+3
\end{matrix}\right).$$ 
$Q$ is also an ED-set, and we have $H_1(Q)=\{3,4,5\}$. So, 
$$Q^2=\left(\begin{matrix}
-4n&0\\
n^2-6n+1&3n+3\\
n^2+2n+1&3n+3
\end{matrix}\right).$$
 Since $H_2(Q)=\emptyset$, $Q^3$ is the empty family. Therefore $w(Q^1)=\left[\begin{matrix} 0&1&2\\
1&2&0\end{matrix}\right]$, $w(Q^2)=\left[\begin{matrix} 1&1\\
1&0\end{matrix}\right]$, and $w(Q^3)=\left[\begin{matrix} 0\\
0\end{matrix}\right]$.

We note that $w(P)=w(Q)$. However, since $w_1(Q)=w_1(P)$ but $w_1(Q^2)<w_1(P^2)$, we have $W(Q)<W(P)$. We have implicitly chosen $D=2$ in this example. As long as $D$ is at least as large as the degree of all polynomial families under consideration, it will not affect whether $W(Q)<W(P)$.

\end{example}

A polynomial family $P= \{p_{ij}\}$ is said to be {\bf
standard} if it is an ED-set and $\deg(p_{1j})=\deg(P)$ for some $1\leq
j\leq l$. We now state Proposition \ref{R} in the case that $P$ is standard.

\begin{proposition}\label{J} Let $(X,\mathcal{B},\mu,(T_1,\ldots,T_l))$ be a freely generated totally
ergodic system and $P = \{p_{ij}\colon 1 \leq i \leq r, 1 \leq j\leq l\}$ be a
standard polynomial family. Then there exists $k\in \bbN$ such that for any $f_1,
\ldots , f_r \in L^{\infty}(\mu)$ with $\nnorm {f_1}_k = 0$, we have
$$\limsup_{N\to\infty} \Big\|\frac{1}{|\Phi_N|}\sum_{u \in \Phi_N} (\prod_{i=1}^r
T_1^{p_{i1}(u)}T_2^{p_{i2}(u)}\ldots T_l^{p_{il}(u)}f_i )\Big\|_{L^2(\mu)}=0$$ for
any F\o lner sequence $\{\Phi_N\}_{N=1}^{\infty}$ in $\bbZ^d$.
\end{proposition}

To prove Proposition \ref{J}, we construct a new polynomial family $Q$ that
controls the above averages, where $W(Q)<W(P)$.  This process is a modified version of the PET-induction process used in \cite{Le4} for a single transformation.

\subsection{Inductive Polynomial Families} We begin by defining that a certain
property holds for {\bf almost all} $v\in \bbZ^d$ if the set of elements for which
the property does not hold is contained in a set of zero density with respect to any F\o lner sequence in $\bbZ^d$. To show a property holds for almost all $v\in\bbZ^d$, we use the fact that a set of zeros of a nontrivial polynomial has zero density with respect to any F\o lner sequence in $\bbZ^d$.

Given any standard polynomial family $P$ with $\deg(P)\geq 2$ where $\deg(p_{11})=\deg(P)$, for each
$(v,w)\in\bbZ^{2d}$ we construct a new family $P_{v,w}$, as follows.  
We first select an appropriate row $i_0$ in $P$, so that $P_{v,w}$ is standard for almost all $(v,w)\in\bbZ^{2d}$ and $W(P_{v,w})<W(P)$.

We split into the following five cases. \label{i_0}
\begin{itemize}
\item {\bf Case 1:} $H_1=\emptyset$ and some $p_{i1}$ is not equivalent to $p_{11}$.

Choose the smallest $i_0$ so that $p_{i_01}$ has minimal degree over all $p_{i1}$ that are not equivalent to $p_{11}$. 
\item {\bf Case 2:} $H_1=\emptyset$, all $p_{i1}$ are equivalent to $p_{11}$, and there is some $i,j$ where $p_{ij}$ is not equivalent to $p_{1j}$ and the degree of either $p_{ij}$ or $p_{1j}$ equals $\deg(P)$.

Choose $i_0$ to be the smallest such $i$ where $p_{ij}$ is not equivalent to $p_{1j}$ and the degree of either $p_{ij}$ or $p_{1j}$ equals $\deg(P)$.
\item {\bf Case 3:} $H_1=\emptyset$, all $p_{i1}$ are equivalent to $p_{11}$ and for all $j$ either $p_{ij}$ is equivalent to $p_{1j}$ for all $i=1\ldots r$.  or $\deg(p_{ij})<\deg(P)$ for all $i=1\ldots r$. 

Choose $i_0=1$.
\item {\bf Case 4:} $H_1\neq \emptyset$, and some $p_{ij_0}$ is not equivalent to $p_{i'j_0}$ for $i,i'\in H_{j_0-1}$.

Choose $i_0$ to be the smallest $i\in H_{j_0-1}$ where $p_{i_0j_0}$ has minimal degree over all $p_{ij_0}$ that are not equivalent to $p_{11}$. 

\item  {\bf Case 5:} $H_1\neq \emptyset$ and all $p_{ij_0}$ are equivalent to $p_{i'j_0}$ for $i,i'\in H_{j_0-1}$.

Choose $i_0=\min H_{j_0-1}$.
\end{itemize}

In our construction, we must treat polynomials in $P$ with degree $1$ differently than those of greater degree. 
For all $(v,w)\in\bbZ^{2d}$, set \label{zij}
$$z_{ij}=\left\{\begin{array}{cc} 
w & \text{ if } \deg(p_{ij})=1\\ v &
\text{ otherwise }\end{array}\right. .$$ For a fixed $(v,w)\in\bbZ^{2d}$,  $p_{ij}(u+z_{ij})$ equals $p_{ij}(u+v)$ or $p_{ij}(u+w)$, depending only on the degree on $p_{ij}$. Thus we view $p_{ij}(u+z_{ij})$ and $p_{ij}(u+w)$ as polynomials in $u$. Given $(v,w)\in\bbZ^{2d}$, we define the new polynomial
family 
\begin{multline*}
\bar{P}_{v,w}=\{p_{ij}(u+z_{ij}),p_{ij}(u+w)\colon i\in I_2,j=1\ldots, l\}\\
\bigcup \{p_{ij}(u+w)\colon i\in I_1,j=1\ldots, l\}.
\end{multline*}

We relabel the family $$\bar{P}_{v,w}=\{q_{v,w,h,j}:1\leq h\leq s,1\leq  j\leq l\}$$ in the
following manner. 
We label each row $$p_{i1}(u+z_{i1}),\ldots,p_{il}(u+z_{il})$$ and
$$p_{i1}(u+w),\ldots,p_{il}(u+w)$$ as $$q_{v,w,h,1}(u),\ldots,q_{v,w,h,l}(u)$$ for
some unique $1\leq h\leq s$ where $p_{1j}(u+z_{1j})=q_{v,w,1,j}$ and
$p_{i_0j}(u+w)=q_{v,w,s,j}(u)$.

Since for each vector $(v,w)$ in $\bbZ^{2d}$, $p_{ij}(u+v),p_{ij}(u+w)$, and
$p_{ij}(u)$ are all equivalent, 
$\bar{P}_{v,w}$ and $P$ have identical column degrees, and $w_j(P)=w_j(P_{v,w})$ for all $1 \leq j\leq l$ and $(v,w)\in
\bbZ^{2d}$. By construction, the first row of $\bar{P}_{v,w}$ also contains a
polynomial of maximal degree and it is easy to check that $\bar{P}_{v,w}$ is
an ED-set for each $(v,w)$ outside a set of zeros of finitely many
polynomials. Hence, $\bar{P}_{v,w}$ is a standard polynomial family for almost all
$(v,w)\in \bbZ^{2d}$.

Next, for each $(v,w)\in\bbZ^{2d}$ we define the new family
$$P_{v,w}=\{q_{v,w,h,j}-q_{v,w,s,j}\colon 1\leq h\leq s-1,1\leq j\leq l\}.$$

\begin{example} For $P$ in our previous example on page \pageref{Example}, case (4) applies and $i_0=2$. It is easy to check that $Q=P_{v,w}$ with $(v,w)=(-1,1)$.
\end{example}

\begin{lemma}\label{Stand}
For each standard polynomial family $P$ where $\deg(P)\geq 2$ and $\deg(p_{11})=\deg(P)$, $P_{v,w}$ is standard for almost all choices of $(v,w)\in\bbZ^{2d}$. Moreover, $C(P_{v,w})\leq C(P)$, and $\deg(P_{v,w})$ equals $\deg(P)$ or $\deg(P)-1$.

\end{lemma}

\begin{proof} Since each entry in $P_{v,w}$ is
constructed by subtracting 2 polynomials from the same column of $\bar{P}_{v,w}$, the maximum degree
in each column  of $P_{v,w}$ cannot increase. Therefore $C(P_{v,w})\leq C(P)$ and \linebreak $\deg(P_{v,w})\leq\deg(P)$. It is easy to check that $P_{v,w}$ is an ED-set whenever $\bar{P}_{v,w}$ is. We now show that the first row in $P_{v,w}$ contains a polynomial of maximal degree. 
 
We split into the five cases used to define $i_0$ on page \pageref{i_0}. In cases 1, 4, and 5,  $p_{i_01}$ is not equivalent to $p_{11}$. When $p_{i_01}$ is not equivalent to $p_{11}$, $$\deg(P_{v,w})\geq\deg(q_{v,w,1,1}-q_{v,w,s,1})=\deg(p_{11})\geq\deg(P_{v,w}).$$ Thus, $\deg(q_{v,w,1,1}-q_{v,w,s,1})=\deg(P_{v,w})$ and the first row in $P_{v,w}$ contains a polynomial of maximal degree. 

 In Case 2,  $p_{i_0j}$ is not equivalent to $p_{1j}$ for some $1\leq j\leq l$ and the degree of either $p_{ij}$ or $p_{1j}$ equals $\deg(P)$. So, $$\deg(P_{v,w})\geq\deg(q_{v,w,1,j}-q_{v,w,s,j})=\deg(P)\geq\deg(P_{v,w}).$$ Thus, $\deg(q_{v,w,1,j}-q_{v,w,s,j})=\deg(P_{v,w})$ and the first row in $P_{v,w}$ contains a polynomial of maximal degree. 

In Case 3, all $p_{i1}$ are equivalent to $p_{11}$, and $i_0=1$. Thus, $$\deg(q_{v,w,1,1}-q_{v,w,s,1})=p_{11}(u+v)-p_{11}(u+w)=\deg(P)-1$$ for almost all  $(v,w)\in\bbZ^{2d}$, since $\deg(p_{11})\geq 2$. Let $j\in \{1,\ldots,l\}$. Then either $p_{ij}$ is equivalent to $p_{1j}$ for all $i=1,\ldots r$ or $\deg(p_{ij})<\deg(P)$ for all $i=1\ldots r$. When $p_{ij}$ is equivalent to $p_{1j}$ , then $$\deg(q_{v,w,h,j}-q_{v,w,s,j})<\deg(p_{1j})\leq \deg(P).$$ When $\deg(p_{ij})<\deg(P)$, $\deg(q_{v,w,h,j}-q_{v,w,s,j})<\deg(P).$ Thus, all polynomials in $P_{v,w}$ have degree less than or equal to $\deg(P)-1$, and  for almost all $(v,w)\in\bbZ^{2d}$, $\deg(q_{v,w,1,1}-q_{v,w,s,1})=\deg(P)-1.$ Therefore the first row in $P_{v,w}$ contains a polynomial of maximal degree.

 In each case, the first row in $P_{v,w}$ contains a polynomial of maximal degree, and $\deg(P_{v,w})$ equals $\deg(P)$ in cases 1,2,4,5 and equals $\deg(P)-1$ in case 3.

\end{proof}

\subsection{Reduction of Weight}
We now show that the above construction leads to a reduction in the weights of our polynomial families.

\begin{proposition}\label{WLoss}
For each $(v,w)\in\bbZ^{2d}$ and  each standard polynomial family $P$ where $\deg(p_{11})=\deg(P)\geq 2$, we have $W(P_{v,w})<W(P)$.
\end{proposition}

\begin{proof}
 We show that $W(P_{v,w})<W(P)$ for each of our five cases used to define $i_0$ on page \pageref{i_0}. 
In cases 1,2, and 3, $p_{i_01}$ has minimal degree over all $p_{i1}$. For all $(v,w)$, the equivalence classes and their degrees in each column remain the same in $\bar{P}_{v,w}$ as in $P$. Thus, $w_1(P)=w_1(\bar{P}_{v,w})$. Column 1 of $P_{v,w}$ is comprised of polynomials $q_{v,w,h,1}-q_{v,w,s,1}$, where $q_{v,w,s,1}$ has minimal degree over all $q_{v,w,h,1}$. We now consider each equivalence class in column 1 of $\bar{P}_{v,w}$ as we pass from $\bar{P}_{v,w}$ to $P_{v,w}$. Each distinct equivalence class in column  $1$ of $\bar{P}_{v,w}$ not containing $q_{v,w,s,1}$, remains a distinct equivalence class of the same degree in column $1$ of $P_{v,w}$. The equivalence class in column one containing $q_{v,w,s,1}$ splits into possibly several equivalence classes of lower degree.  Thus, $w_1(P_{v,w})<w_1(P)$, and hence $W(P_{v,w})<W(P)$.

For cases 4 and 5, we show that $w_1((P_{v,w})^{j_0})<w_1(P^{j_0})$, and \linebreak$w_1((P_{v,w})^{a})<w_1(P^{a})$ for all $a<j_0$. The polynomials in the first column of $P^a$ are precisely those entries in the $a^{th}$ column of $P$ only from those rows whose entries are zero in columns $1,\ldots, a-1$. Thus, $w_1(P^a)$ counts the equivalence classes of polynomials from only those rows of column $a$ in $P$ whose entries are zero in columns $1,\ldots,a-1$.

Suppose $1\leq a\leq j_0$. If the $h^{th}$ row of $\bar{P}_{v,w}$ has zeros in columns $1,\ldots,a-1$, then $q_{v,w,h,a}= p_{ia}(u+v) \text{ or } p_{ia}(u+w)$ where $p_{ia}(u)$ is a polynomial in $P$ with $i\in H_{a-1}$. Moreover, for each $i\in H_{a-1}$,  there is some row $h$ of $\bar{P}_{v,w}$ with zeros in columns $1,\ldots,a-1$ and $q_{v,w,h,a}=p_{ia}(u+w)$. Thus, the equivalence classes in $\bar{P}_{v,w}$ from only those rows of column $a$ whose entries are zero in columns $1,\ldots,a-1$  are the same as the equivalence classes in $P$ from only those rows of column $a$ whose entries are zero in columns $1,\ldots,a-1$. Thus, $w_1(\bar{P}^a_{v,w})=w_1(P^a)$.

Since $i_0\in H_{j_0-1}$, $q_{v,w,s,j}=0$ for all $j=1,\ldots, j_0-1$. So, for all $j=1,\ldots,j_0-1$, $q_{v,w,h,j}-q_{v,w,s,j}=q_{v,w,h,j}$. Thus the rows in $\bar{P}_{v,w}$ (except the last) with zeros in columns $1,\ldots,a-1$, are the same as the rows in $P_{v,w}$ with zeros in columns $1,\ldots,a-1$. 

When $a<j_0$, we have $q_{v,w,s,j}=0$, for all $j=1,\ldots, a$. So the equivalence classes and their degrees in only those rows of column $a$ whose entries are zero in columns $1,\ldots,a-1$ are the same for both $\bar{P}_{v,w}$ and $P_{v,w}$. Therefore, $w_1(P^{a}_{v,w})=w_1(\bar{P}_{v,w}^a)=w_1(P^a)$.

When $a=j_0$, we have $q_{v,w,s,a}\neq 0$. However, $q_{v,w,s,a}$ has minimal degree over all $q_{v,w,h,a}$ where $q_{v,w,h,a}=0$ for all $j=1,\ldots, a-1$. As before, each distinct equivalence class of such polynomials in column  $a$ of $\bar{P}_{v,w}$ not containing $q_{v,w,s,a}$, remains a distinct equivalence class of the same degree in column $a$ of $P_{v,w}$. The equivalence class in column $a$ containing $q_{v,w,s,a}$ splits into possibly several equivalence classes of lower degree. Therefore, $w_1(P^{a}_{v,w})<w_1(P^a)$.
Since,  $w_1(P^{a}_{v,w})=w_1(P^a)$ for $a=1,\ldots, j_0-1$ and $w_1(P^{j_0}_{v,w})=w_1(P^{j_0})$, $W(P_{v,w})<W(P)$.
\end{proof}

\subsection{PET-Induction}
\begin{proof}[Proof of Proposition \ref{J}] 
Let $P=\{p_{ij}\colon 1 \leq i \leq r, 1 \leq j\leq l\}$ be a standard polynomial family.
For polynomial families of degree 1, the result is
given by Proposition \ref{G}. Suppose $\deg(P)\geq 2$. Since $P$ is standard, by relabeling the transformations, we may assume that $deg(p_{11})=\deg(P)$. There are only finitely many column degrees $C(Q)<C(P)$ and weights $W(Q)<Q(P)$ that correspond to families $Q=\{q_{ij}\colon 1 \leq i \leq s, 1 \leq j\leq l\}$ where $1\leq s\leq 2r$ and $C(Q)\leq C(P)$. Thus, we state our PET-induction hypothesis as follows. We assume that for all $1\leq s\leq 2r$ there exists $k\in\bbN$ such that for all standard polynomial families $Q=\{q_{ij}\colon 1 \leq i \leq s, 1 \leq j\leq
l\}$ where $C(Q)<C(P)$,  or where $C(Q)\leq C(P)$, $\deg(q_{11})=\deg(Q)$, and $W(Q)<W(P)$, we have $$\limsup_{N\to\infty} \Big\|\frac{1}{|\Phi_N|}\sum_{u \in \Phi_N} 
(\prod_{i=1}^{s} T_1^{q_{i1}(u)}\ldots T_l^{q_{il}(u)}b_i )\Big\|_{L^2(\mu)}=0,$$ for any $b_1, \ldots , b_r \in L^{\infty}(\mu)$ with $\nnorm {b_1}_k = 0$, and for each F\o lner sequence $\{\Phi_N\}_{N=1}^{\infty}$ in $\bbZ^d$.

Now let $f_1,\ldots,f_r \in L^{\infty}(\mu)$ where  $\nnorm{f_1}_k=0$, and let $\{\Phi_N\}_{N=1}^{\infty}$ be
a F\o lner sequence in $\bbZ^d$. Without loss of generality we may assume that
$\norm{f_i}_{L^{\infty}(\mu)}\leq 1$ for all $1\leq i\leq r$. By replacing each $f_i$ with $T_1^{c_1}\ldots T_l^{c_l}f_i$ for some $c_1,\ldots, c_l\in\bbZ$, we may assume that each $p_{ij}$ has zero constant term. In particular, each polynomial in $P$ whose degree is $1$ is linear.

By Lemma \ref{F} and the Cauchy-Schwartz inequality we have for any finite set $F\subset \bbZ^d$,
\begin{multline*}
\limsup_{N\to\infty} \Big\|\frac{1}{|\Phi_N|}\sum_{u \in \Phi_N} 
(\prod_{i=1}^{r} T_1^{p_{i1}(u)}\ldots T_l^{p_{il}(u)}f_i )\Big\|_{L^2(\mu)}^2 \\
\leq \limsup_{N\to\infty} \frac{1}{|F|^2}\sum_{v,w\in F}\frac{1}{|\Phi_N|}
\sum_{u\in \Phi_N}\int_X\prod_{i=1}^{r} T_1^{p_{i1}(u+v)}
\ldots \\
T_l^{p_{il}(u+v)}f_i \cdot \prod_{i=1}^{r} T_1^{p_{i1}(u+w)}\ldots
T_l^{p_{il}(u+w)}f_i d\mu \\
\leq \limsup_{N\to\infty} \frac{1}{|F|^2}\sum_{v,w\in
F}\frac{1}{|\Phi_N|}\sum_{u\in \Phi_N}\int_X\prod_{h=1}^{s}
T_1^{q_{v,w,h,1}(u)}\\ 
\ldots T_l^{q_{v,w,h,l}(u)}b_{v,w,h}d\mu \\
\leq \frac{1}{|F|^2}\sum_{v,w\in
F}\limsup_{N\to\infty}\Big\|\frac{1}{|\Phi_N|}\sum_{u\in \Phi_N}
\prod_{h=1}^{s-1}T_1^{(q_{v,w,h,1}-q_{v,w,s,1})(u)} \\
 \ldots
T_l^{(q_{v,w,h,l}-q_{v,w,s,l})(u)}b_{v,w,h}\Big\|_{L^2(\mu)}
\end{multline*}
for $(v,w) \in \bbZ^{2d}$, where the $b_{v,w,h}$ represent any of
the following bounded functions: 
\begin{itemize}
\item $T_1^{p_{i1}(v-z_{i1})}\ldots T_l^{p_{il}(v-z_{il})}f_i$ for $i\in I_2$,
\item $f_i\cdot T_1^{p_{i1}(v)-p_{i1}(w)}\ldots
T_l^{p_{il}(v)-p_{il}(w)}f_i$ for $i\in I_1$.
\end{itemize}
Since $P$ has degree of at least $2$, $1\in I_2$ and $b_{v,w,1}=T_1^{t_1}\ldots T_l^{t_l}f_1$ for some $t_1,\ldots,t_l\in\bbZ$. Thus, for all $k\in\bbN$ and all $(v,w)\in\bbZ^{2d}$, 
$$\nnorm{b_{v,w,1}}_k=\nnorm{f_1}_k.$$
However, $P_{v,w}=\{q_{v,w,h,j}-q_{v,w,s,j}\colon 1\leq h\leq s-1,1\leq j\leq l\},$ is a standard polynomial family where $1\leq s-1\leq 2r$ and $W(P_{v,w})<W(P)$ for almost all $(v,w)\in\bbZ^{2d}$. We note that whenever $\deg(q_{v,w,1,1}-q_{v,w,s,1})<\deg(P_{v,w})$, $C(P_{v,w})<C(P)$. By the PET-induction hypothesis, for almost all choices  of $(v,w)\in\bbZ^{2d}$, we have
$$\limsup_{N\to\infty}\Big\|\frac{1}{|\Phi_N|}\sum_{u\in \Phi_N}
\prod_{h=1}^{s-1}T_1^{(q_{v,w,h,1}-q_{v,w,s,1})(u)} \\
 \ldots
T_l^{(q_{v,w,h,l}-q_{v,w,s,l})(u)}b_{v,w,h}\Big\|_{L^2(\mu)}=0.$$ For all other choices of $(v,w)\in\bbZ^{2d}$, the above average is bounded above by $1$. Therefore,

\begin{multline*}
\limsup_{N\to\infty} \Big\|\frac{1}{|\Phi_N|}\sum_{u \in \Phi_N} 
(\prod_{i=1}^{r} T_1^{p_{i1}(u)}\ldots T_l^{p_{il}(u)}f_i )\Big\|_{L^2(\mu)}^2\\
\leq \inf_{F}\frac{1}{|F|^2}\sum_{v,w\in
F}\limsup_{N\to\infty}\Big\|\frac{1}{|\Phi_N|}\sum_{u\in \Phi_N}
\prod_{h=1}^{s-1}T_1^{(q_{v,w,h,1}-q_{v,w,s,1})(u)} \\
 \ldots
T_l^{(q_{v,w,h,l}-q_{v,w,s,l})(u)}b_{v,w,h}\Big\|_{L^2(\mu)}=0
\end{multline*}

where the infimum is taken over all finite subsets of $\bbZ^{d}$.
\end{proof}

\subsection{Reduction to the standard case}
\begin{proof}[Proof of Proposition \ref{R}]
We now reduce the general case to one involving standard systems. Let $P = \{p_{ij}
\colon 1 \leq i \leq r, 1 \leq j\leq l\}$ be a (nonstandard) ED-set of polynomials of degree less than $b$, let $f_1,\ldots,f_r\in
L^{\infty}(\mu)$, and let $\{\Phi_N\}_{N=1}^{\infty}$ be a F\o lner sequence in
$\bbZ^d$. Once again, we assume that each polynomial in $P$ has zero constant term. In otherwords, $p_{ij}({\bf 0})=0$ for each polynomial $p_{ij}$ in $P$, where ${\bf 0}$ is the zero vector in $\bbZ^{d}$.
Thus, we have $p_{ij}(u+v)=p_{ij}(u+z_{ij})+p_{ij}(v-z_{ij})$ for each polynomial in $P$, where $z_{ij}$ is defined as on page \pageref{zij}.
By Lemma \ref{F}, there exists a F\o lner sequence
$\{\Theta_N\}_{N=1}^{\infty}$ in $\bbZ^{3d}$ such that 
\begin{multline*}
\limsup_{N\to\infty} \Big\|\frac{1}{|\Phi_N|}\sum_{u \in \Phi_N} 
(\prod_{i=1}^r T_1^{p_{i1}(u)}\ldots T_l^{p_{il}(u)}f_i )\Big\|_{L^2(\mu)}^2 \\
\leq \limsup_{M\to\infty}\frac{1}{\Theta_M}\sum_{(u,v,w)\in\Theta_M}\int_X
\prod_{i=1}^r T_1^{p_{i1}(u+v)+q(u)}\\
\ldots T_l^{p_{il}(u+v)}f_i \prod_{i=1}^r T_1^{p_{i1}(u+w)+q(u)}\ldots
T_l^{p_{il}(u+w)}f_i d\mu \\
\leq \limsup_{M\to\infty}\Big\|\frac{1}{\Theta_M}\sum_{(u,v,w)\in\Theta_M}
\prod_{i=1}^r T_1^{p_{i1}(u+z_{i1})+q(u)}\ldots
T_l^{p_{il}(u+z_{il})}(T_1^{p_{i1}(v-z_{i1})}\\
\ldots T_l^{p_{il}(v-z_{il})}f_i) \prod_{i=1}^r T_1^{p_{i1}(u+w)+q(u)}\ldots
T_l^{p_{il}(u+w)}f_i\Big\|_{L^2(\mu)} \\
\end{multline*}
where $q:\bbZ^{3d}\to\bbZ$ is any polynomial of degree $b$. Whether $z_{ij}$ equals $v$ or $w$ is determined only by the degree of $p_{ij}$, so each polynomial below is really only a polynomial in $u,v,w$. Thus the set 
$$\{p_{i1}(u+z_{i1})+q(u),p_{i1}(u+w)+q(u),p_{ij}(u+z_{ij}),p_{ij}(u+w)\colon 1\leq i\leq
r, 2\leq j\leq l\}$$ 
of polynomials from $\bbZ^{3d}\to\bbZ$ is a standard family of degree $b$. Thus there
exists $k\in \bbN$ (that depends only on the original polynomial family $P$) such that

\begin{multline*}
\limsup_{M\to\infty}\Big\|\frac{1}{\Theta_M}\sum_{(u,v,w)\in\Theta_M} \prod_{i=1}^r
T_1^{p_{i1}(u+z_{i1})+q(u)}\ldots T_l^{p_{il}(u+z_{il})}(T_1^{p_{i1}(v-z_{i1})}\\
\ldots T_l^{p_{il}(v-z_{il})}f_i) \prod_{i=1}^r T_1^{p_{i1}(u+w)+q(u)}\ldots
T_l^{p_{il}(u+w)}f_i\Big\|_{L^2(\mu)}=0.
\end{multline*}
\end{proof}

\end{document}